\documentclass[12pt]{article}
\usepackage{amssymb}

\def\no{\nonumber}
\def\be{\begin{equation}}
\def\ee{\end{equation}}
\def\ba{\begin{eqnarray}}
\def\ea{\end{eqnarray}}

\def\btu{\bigtriangleup}
\def\btd{\bigtriangledown}
\def\e1{\epsilon}
\def\o1{\Omega}
\def\a1{\alpha}
\def\om{\omega _{\e1}}
\def\l1{\Lambda}
\def\xl{\lambda_{\e1}}
\def\sl{\lambda_{\e1}^*}
\def\v1{\varphi}
\def\xv{\varphi_{\e1}}
\def\0v{\varphi_{0}}
\def\d1{\delta}
\def\p1{\partial}
\def\ph{\phi _{\e1}}
\def\f1{\frac}
\def\t1{\theta}
\def\s1{\sqrt{\e1}}
\def\b1{\beta (r,\theta )}
\def\bar{\overline}
\def\bs{\begin{eqnarray*}}
\def\es{\end{eqnarray*}}
\begin{document}
\title{The differential equation $\btu
 u = 8\pi - 8\pi he^u$ on a
compact Riemann surface\thanks{{\bf dg-ga/9710005} }}
\author{Weiyue Ding\thanks{Institute of Mathematics,
Academia Sinica, Beijing 100080, P. R. China.
Partially supported by CNNSF.}\\
\and
J\"{u}rgen Jost\thanks{Max-Planck-Institute for Mathematics in the Sciences,
Inselstr. 22-26, 04103 Leipzig, Germany.}\\
\and
Jiayu Li\thanks{Institute of Mathematics,
Academia Sinica, Beijing 100080, P. R. China.}\\
\and
Guofang Wang\thanks{Institute of Systems Science,
Academia Sinica, Beijing 100080, P. R. China.}}
\date {June 1997}

\maketitle

\begin{abstract}

Let $M$ be a compact Riemann surface, $h(x)$ a positive smooth function on $M$.
In this paper, we  consider the functional
$$
J(u)=\frac{1}{2}\int _{M}|
\btd 
 u|^2 + 8\pi \int _{M}u
      -8\pi \log\int _{M}he^{u}.
$$
We give a sufficient condition under which $J$ achieves its minimum.
\end{abstract}

\newtheorem{theorem}{Theorem}[section]
\newtheorem{lemma}[theorem]{Lemma}
\newtheorem{corollary}[theorem]{Corollary}
\newtheorem{remark}[theorem]{Remark}
\newtheorem{definition}[theorem]{Definition}
\newtheorem{proposition}[theorem]{Proposition}

\section{Introduction and main result}

Let $(M,ds^2)$ be a compact Riemann surface, $h(x)$ a smooth function on $M$.
For simplicity, we assume in this paper that the volume of $M$ equals
1. Twenty years ago, Kazdan and Warner ([KW]) asked, under what
kind of conditions on $h$, the equation
\be
\btu u = 8\pi - 8\pi he^u
\ee
has a solution. An obvious necessary condition is that $\max h >0$.

If $M$ is the standard sphere, the problem is called ``Nirenberg problem''.
The geometric significance of this problem is that if $g$ denotes a metric
of constant curvature $4\pi$ on $S^2$, then the metric $e^{u}g$ has
curvature equal to $h$.
This problem has been studied by Moser ([M1], [M2]), Kazdan-Warner ([KW]),
Hong ([H]),
Chen-Ding ([CD1], [CD2]), Chang-Yang ([CY1], [CY2]), Chang-Liu ([CL]),
and others.

For a compact Riemann surface other than $S^2$ or $\mathbb{RP}^2,$ the preceding
interpretation is no longer possible as such a surface does not carry a
background
metric of constant positive curvature. However, the differential equation (1)
also arises in the so-called Chern-Simons-Higgs theory. This is a classical
field theory that is defined on (2+1) dimensional Minkowski space and believed
to be relevant in high temperature superconductivity and in other areas of
theoretical physics. Hong-Kim-Pac [HKP] and Jackiw-Weinberger [JW] observed
that for a special choice of the Higgs potential, a sixth order polynomial,
stationary vortex solutions satisfy certain first order selfduality equations.
On a compact torus, these equations have been studied by Caffarelli-Yang [CaY]
and Tarantello [T]. In particular, Tarantello showed that one may find a
certain type of solution that corresponds to a symmetric vacuum. In the case
of only one vortex $p$ of multiplicity 1 she found that asymptotically,
as the coupling parameter in the theory tends to zero, one obtains a solution
of
$$
\btu u(x) = 4\pi - 4\pi \frac{e^{-G(x,p)+u(x)}}{\int_{M}e^{-G(y,p)+u(y)}dy},
\quad \quad
\int_{M}u = 0, \quad \quad u \in L^2_1 (M),
$$
where $G(x,p)$ is the Green function defined below in equation (2).
This result was shown with the help of the Moser-Trudinger inequality.
For $N$ vortices (counted with multiplicity), $4\pi$ in the preceding
equation has to be replaced by $4\pi N,$ and already for $N=2$, the
situation becomes more difficult, as the factor $8\pi $ represents the
limiting case in the Moser-Trudinger inequality. We shall discuss this in
detail in [DJLW].

Thus, the Kazdan-Warner problem becomes relevant for an area quite different
from problems of prescribed Gauss curvature. Therefore, we shall address the
problem of finding solutions (1) here on a general compact Riemann surface.
We shall pursue a variational approach. Namely, we shall try to minimize the
functional
$$
J(u)=\frac{1}{2}\int _{M}|\btd u|^2 + 8\pi \int _{M}u
      -8\pi \log\int _{M}he^{u}.
$$

We shall first show that the functional has a lower bound. This generalizes the
Moser-Trudinger ([M1]) inequality to the case where $M$ is an arbitrary
compact Riemann surface.

\begin{theorem}
Let $(M,ds^2)$ be a compact Riemann surface. For any $u\in L_1^2(M)$ with
$\int _M u=0$ one has

$$
\int _M e^u \leq C_M e^{\frac{1}{16\pi}\|\btd u\|_2^2},
$$
where $C_M$ is a positive constant depending only on $(M, ds^2)$.
\end{theorem}

To show that $J$ is bounded from below, we consider
$$
J_{\e1 }(u)=\frac{1}{2}\int _{M}|\btd u|^2 + (8\pi -\e1 )\int _{M}u
-(8\pi -\e1 ) \log\int _{M}he^{u},
$$
where $\e1 >0$. It is not hard to verify that $J_{\e1 }$ achieves its
minimum at some $u_{\e1}.$

There are two possibilities: If a subsequence of the sequence of minimizers
$u_{\e1}$ converges to some $u_0$ for $\e1 \rightarrow 0$, then $u_0$ minimizes
$J$. In order to show this convergence, it suffices to establish estimates
for  $u_{\e1}$ in the Sobolev space $L^2_1 (M)$ that do not depend on
$\e1.$ If such estimates do not hold, then the sequence $u_{\e1}$ blows up,
and after subtracting mean values, the $u_{\e1}$ converges to some Green
function $G(x,p)$ satisfying

\be
\cases{
     \btu G = 8\pi -8\pi\delta _p,\cr
     \int _M G=0.}
\ee

In a normal coordinate system around $p$ we assume that
\ba
G(x,p) &=&-4\log r + A(p) + b_1x_1 + b_2x_2\no \\
&&+c_1x_1^2 + 2c_2x_1x_2 + c_3x_2^2 + O(r^3),
\ea
where $r(x)=dist (x,p)$.

One should note that (2) is not conformally invariant, but depends on the
metric $ds^2$ on $M$. Therefore, also the constants in the expansion (3)
will depend on that metric. If the metric is homogeneous as on the standard
sphere or on a flat torus, $b_1 = b_2 = 0.$ For a more detailed discussion
of the leading term $A(p)$ - which does not depend on $p$ in the homogeneous
case - on flat tori see section 4.

More precisely, in this step we show that, if the minimizing sequence
$u_\e1$ of $J_{\e1 }$ blows up,
\be
\inf _{u\in L^2_1(M)} J(u)\geq -8\pi -8\pi\log\pi - 4\pi
( \max_{p\in M} (A(p)+2\log h(p))).
\ee

In other words, if (4) does not hold, then no blow-up is possible,
and we get convergence of the $u_{\e1}$ to a minimizer $u_0$ of $J.$

Inequality (4) and the results that have been obtained for the Nirenberg
problem ([CD1], [CY1], [CY2], [CL]) indicate that it will depend on the
asymptotic expansion of $h$ near a potential blow-up point whether a
blow-up is possible. In this sense, we shall obtain the following result.

\begin{theorem}
Let $(M,ds^2)$ be a compact Riemann surface, let $K(x)$ be its Gauss curvature.
Let $h(x)$ be a positive smooth function on $M$.
Suppose that $A(p)+2\log h(p)$ achieves its maximum at $p_0$.
Let $b_1(p_0)$ and $b_2(p_0)$ be the constants in the expression (3),
and write $\btd h(p_0) = (k_1(p_0),k_2(p_0))$ in the normal coordinate system.
If
\bs
\lefteqn{\btu h (p_0) + 2(b_1(p_0)k_1(p_0)+b_2(p_0)k_2(p_0))}\\
& > & -(8\pi + (b_1^2(p_0)+b_2^2(p_0)) - 2K(p_0))h(p_0)
\es
the minimum of the functional $J$ can be obtained, and consequently
the equation (1) has a smooth solution.
\end{theorem}

\begin{remark}
The inequality in Theorem 1.2 is implied by the following one
\bs
\btu \log h(p_0)=\f1{\btu h (p_0)}{h(p_0)} - \frac{| \btd h (p_0)|^2}{h^2
(p_0)}
>
- (8\pi - 2K(p_0)).
\es
\end{remark}

In the second step, we shall construct a blowing up sequence $\ph $ with
the property that
$$
J(\ph ) < -8\pi -8\pi\log\pi - 4\pi
( \max_{p\in M} (A(p)+2\log h(p)))
$$
for sufficiently small $\e1 >0$,
assuming that $h$ satisfies the hypotheses in Theorem 1.2. This contradicts
(4), and Theorem 1.2 will follow.

 Our methods are closely related to those used by Schoen ([Sc]) in his
 solution of the Yamabe problem and by Escobar-Schoen ([E-S]) for finding
 conformal metrics with prescribed curvatures in higher dimensions. However,
 our analysis is more delicate. In their work, they need only to compare the
 minimum of the corresponding functional to the minimum on the standard
 sphere. That is because their problems are conformally invariant. In
 our case, we have to compute the limit functional value of a
 blowing up minimizing sequence very carefully, and it turns out that
 the limit is not unique, it depends on the geometry of the surface
 (Theorem 1.2). On the other hand, while in their work to establish the
 existence result they need only the constant term in the expansion of
 the Green function of the conformal Laplacian to be positive (the positive
 mass theorem), in our case we need to consider a higher order term in
 the expansion of the usual Green function.

{\bf Acknowledgement:}\\
Our research was carried out at the Max-Planck-Institute for
Mathematics in the Sciences in Leipzig. The first author thanks the
Max-Planck-Institute for the hospitality and good working conditions.
The third author was supported
by a fellowship of the Humboldt foundation, whereas the fourth author
was supported by the DFG through the Leibniz award of the second author.

\section{The lower bound}

In this section, we shall show that the functional $J(u)$ is bounded
from below, and consequently, we shall prove the Moser-Trudinger inequality.

We shall consider the minimum of the functional $J$ in the space
$H_1 =\{~u\in L^2_1(M)~|~\int _M he^u =1~\}$.

\begin{proposition}
Let $M$ be a compact Riemann surface.
Let $h(x)$ be a positive smooth function on $M$.
Then there exists a positive constant
$C$ depending only on $M$ and $h$ such that
$$
\inf_{u\in H_1} J (u)\geq -C.
$$
\end{proposition}

The following lemma will yield our proposition.

\begin{lemma}
There exists a positive constant
$C$ depending only on $M$ and $h$, but not on $\e1$, such that
$$
\inf_{u\in H_1} J_{\e1}(u)\geq -C.
$$
\end{lemma}
Set
$$
\l1 _{\e1 } = \inf_{u\in H_1} J_{\e1}(u).
$$
Using Aubin's inequality (see [A])
one obtains $u_{\e1}\in H_1$ satisfying
$$J_{\e1}(u_{\e1})=\l1 _{\e1 }$$
and
\be
\btu u_{\e1}=(8\pi -\e1 )-(8\pi -\e1 )he^{u_{\e1}}.
\ee

If $u_{\e1}\to u_0$ in $L^2_1(M)$ as $\e1 \to 0$, the lemma follows.
And Theorem 1.2 also follows. Therefore we shall assume in the sequel that
$u_{\e1}$ does not converge in $L^2_1(M)$. However we have

\begin{lemma}
For any $1<q<2$, $\|\btd u_{\e1} \|_q\leq C_q $.
\end{lemma}

Proof:
Let $q'=\frac{q}{q-1}>2$. Then
$$
\|\btd u_{\e1} \|_q\leq \sup \{~|\int _M \btd u_{\e1}\cdot \btd \varphi |~
|~\v1\in L_1^{q'}(M), \int_M\v1 =0, \|\v1 \|_{L_1^{q'}(M)}=1~\}.
$$

By the Sobolev embedding theorem we have
$$
\|\v1 \|_{L^{\infty}(M)}\leq C.
$$
It is clear that
$$
|\int _M \btd u_{\e1}\cdot \btd \varphi |
=|\int _M \btu u_{\e1}\varphi |\leq C.
$$
This proves the lemma.

\hfill $\square$

Let $\bar{u}_{\e1}=\int _M u_{\e1}$.
We set $\xl =\max _{x\in M}u_{\e1}(x)$, assume that
$u_{\e1}(x_{\e1})=\xl$ and that $x_{\e1}\to p$.
We shall show

\begin{lemma}
$\xl \to \infty$ as $\e1\to 0$.
\end{lemma}

Proof:
If $\xl$ did not tend to $\infty$, $ e^{u_{\e1}} $ would be bounded above
(At least there would exist a subsequence $u_{\e1 _{k}}$ such that
$ e^{u_{\e1 _{k}}} $ is bounded. For simplicity, in this paper we do not
distinguish this point.). We set $v_{\e1}=u_{\e1}-\bar{u}_{\e1}$.
By Lemma 2.3 we have $\|v_{\e1}\|_{p}\leq C_p$ for any $p>1$. Since
$|\btu v_{\e1}|\leq C$, by the elliptic estimates we can see that
$v_{\e1}$ is bounded in $C^k(M)$. So, if $\bar{u}_{\e1}$ is bounded,
then $\|u_{\e1}\|_{L^{\infty}(M)}\leq C$, which contradicts the assumption
that $u_{\e1}$ blows up. But if $\bar{u}_{\e1}\to -\infty$, then $v_{\e1}$
converges to a smooth solution of the equation $\btu v = 8\pi$, which
is impossible. This proves the lemma.

\hfill $\square$

We choose a local
normal coordinate system around $p$. Let $(\sl)^2=e^{\xl}$, and
$$
\xv (x)= u_{\e1}(x_{\e1} + \frac{x}{\sl})-\xl.
$$
We shall show
\begin{lemma}
$(i)$ For any $\o1\subset\subset M\setminus \{p\}$, we have
$\int _{\o1}he^{u_{\e1}}\to 0$ as $\e1 \to 0$. And $\bar{u}_{\e1}\to
-\infty$ as $\e1\to 0$.
$(ii)$ For any $\o1\subset\subset R^2$,
we have $\xv (x)\to \v1 _0(x)$ in $C^{\infty}(\o1 )$ as $\e1\to 0$,
where $\v1 _0(x)=-2\log (1+\pi h(p)|x|^2)$.
\end{lemma}

Proof:
For any $R>0$, we have
\bs
\btu\xv &=& \frac{1}{(\sl )^2}(8\pi -\e1 )-(8\pi -\e1 )
\frac{h(x_{\e1}+\frac{x}{\sl})}{(\sl )^2}
e^{\xv (x)+\xl }\\
&=&\frac{1}{(\sl )^2}(8\pi -\e1 )-(8\pi -\e1 )
h(x_{\e1}+\frac{x}{\sl})e^{\xv (x)}
\es
in $B_R (0)\subset R^2$ for $\e1 >0$ sufficiently small.

We consider the equation
$$
\left\{ \begin{array}{ll}
 \btu \xv ^1=  \frac{1}{(\sl )^2}(8\pi -\e1 )-(8\pi -\e1 )
h(x_{\e1}+\frac{x}{\sl})e^{\xv (x)}
 & \quad x\in B_{R}(0),\\
 \xv ^1|_{\partial B_{R}(0)}=0\\
\end{array}
\right.
$$
Let $\xv ^2=\xv -\xv ^1$. Then $\btu \xv ^2=0$ in $B_{R}(0)$.
The elliptic estimates together with $L_2^p(B_{R}(0))\subset
C(\bar{B_{R}(0)})$ give
$
\sup _{B_{R}(0)}|\xv^1|\leq  C.
$
So
$
\sup_{B_{R}(0)}\xv ^2 \leq  C.
$
The Harnack inequality yields that
$
\sup_{B_{\f1{R}{2}}(0)}|\xv ^2 | \leq  C,
$
because $\xv ^2 (0)$ is bounded.
Therefore
$
\sup_{B_{\frac{R}{2}}(p)}|\xv | \leq  C.
$

By the elliptic estimates, we
can show that $\xv (x)\to \v1 _0(x)$ in $C^{\infty}(B_{\frac{R}{4}}(0) )$.
As
$h(x_{\e1}+\frac{x}{\sl})\to h(p)$ in $C(B_R(0))$
we can see that $\0v$ satisfies
$$
\btu_0 \0v (x) = -8\pi h(p)e^{\0v},
$$
$$\0v (0)=0,$$
and
$$
\int _{R^2}h(p)e^{\0v}\leq 1,
$$
where $\btu_0$ is the Laplace operator on $R^2$.

However Ding's lemma ([D], c.f. [CL2] Lemma 1.1) yields that
$$
\int _{R^2}h(p)e^{\0v}=1.
$$

Since $\int _{M}he^{u_{\e1}}=1$ we can see that,
for any $\o1\subset\subset M\setminus \{p\}$, we have
$\int _{\o1}he^{u_{\e1}}\to 0$ as $\e1 \to 0$.
By Jessen's inequality we have $\bar{u}_{\e1}\to
-\infty$ as $\e1\to 0$. The uniqueness theorem in [CL2] implies that
$$
\v1 _0(x)=-2\log (1+\pi h(p)|x|^2).
$$
This proves the lemma.

\hfill $\square$

We also need the following lemma.
\begin{lemma}
We have $u_{\e1}-\bar{u}_{\e1}
\to G(x,p)$ weakly in $L^q_1(M)$ $(1<q<2)$ as $\e1 \to 0$, where
$G$ is the Green function satisfying $(2)$, $p\in M$. Furthermore
$u_{\e1}-\bar{u}_{\e1}\to G(x,p)$ in $C^{\infty}(\o1)$ for any
$\o1 \subset\subset M\setminus \{p\}$.
\end{lemma}

In order to show the lemma, we shall use a theorem proved
by Brezis and Merle ([BM], Theorem 1), formulated as Lemma 2.7 below.

Let $\o1 \subset M$ be a domain. Suppose that $u$ is a
solution of the equation
$$
\left\{ \begin{array}{ll}
 \btu u = f(x)\\
 u|_{\partial\o1 }=0
\end{array}
\right.
$$
in $\o1$ with $\|f\|_{L^1(\o1)}<\infty$.

\begin{lemma}
For any $0<\d1 <4\pi$, we have
$$
\int _{\o1}\exp\{\frac{(4\pi -\d1 )|u(x)|}{\|f\|_{L^1(\o1 )}}\}
\leq C_{\d1 },
$$
where $C_{\delta}>0$ is independent of $\|f\|_{L^1(\o1 )}$.
\end{lemma}

\hfill $\square$

Using Lemma 2.7 we can show

\begin{lemma}
Suppose that $\o1\subset M$ is a domain. If
$$
\int _{\o1}he^{u_{\e1}}\leq (\frac{1}{2}-\d1 )
$$
for some $0<\d1 <\frac{1}{2}$, then
$$
\|u_{\e1}-\bar{u}_{\e1}\|_{L^{\infty}(\o1 _0 )}\leq C(\o1 _0,\o1 )
$$
for any $\o1 _0\subset\subset \o1$.
\end{lemma}

Proof:
Assume that $u_{\e1}^1$ is a solution of the equation
$$
\left\{ \begin{array}{ll}
 \btu  u_{\e1}^1 = -(8\pi -\e1 )he^{u_{\e1}} ,\\
 u_{\e1}^1|_{\partial\o1 }=0 .
\end{array}
\right.
$$

Set $u_{\e1}^2=u_{\e1}-u_{\e1}^1-\bar{u}_{\e1}$, then
$\btu u_{\e1}^2=(8\pi -\e1 )$ in $\o1$. Harnack's inequality
yields that
\bs
\|u_{\e1}^2\|_{L^{\infty}(\o1 _1 )}
& \leq & C(\|u_{\e1}^2\|_{L^1(\o1 )})\\
& \leq & C(\|u_{\e1}-\bar{u}_{\e1}\|_{L^1(\o1 )}
+ \|u_{\e1}^1 \|_{L^1(\o1 )})\\
& \leq & C(\|\btd u_{\e1}\|_{L^q(M )}
+ \|u_{\e1}^1\|_{L^1(\o1 )})
\es
whenever $\o1 _0\subset\subset\o1 _1\subset\subset\o1$.

By Lemma 2.7, one can see that $e^{|u_{\e1}^1 |}$ is bounded
in $L^p(\o1 )$ for some $p>1$, which yields that
$$
\|u_{\e1}^1\|_{L^1(\o1 )}\leq C.
$$
We therefore have
$$
\|u_{\e1}^2\|_{L^{\infty}(\o1 _1 )}\leq C.
$$

Note that
\bs
\int _{\o1 _1}e^{pu_{\e1}} & = &\int _{\o1 _1}e^{p\bar{u}_{\e1}}e^{pu_{\e1}^2}
e^{pu_{\e1}^1}\\
& & \leq C \int _{\o1 _1}e^{p|u_{\e1}^1|}\\
& & \leq C.
\es
By the standard elliptic estimates, we can obtain
$$
\|u_{\e1}^1\|_{L^{\infty}(\o1 _0 )}\leq C.
$$
This proves the lemma.

\hfill $\square$

Now we turn to the proof of Lemma 2.6.

{\it Proof of Lemma 2.6}: By Lemma 2.5 we can see that
$(8\pi -\e1 )he^{u_{\e1}}$
converges to $8\pi\d1 _p$ in the sense of measures as $\e1\to 0$.

Therefore $u_{\e1}-\bar{u}_{\e1}\to G(x,p)$ weakly in $L_1^q(M)$ for
any $1<q<2$, where $G$ is the Green function satisfying (2), because
$G$ is the only solution of (2) in $L_1^q(M)$.

Lemma 2.5 and Lemma 2.8 yield that for any
$\o1\subset\subset M\setminus \{p\}$,
\be
\|u_{\e1}-\bar{u}_{\e1}\|_{L^{\infty}(\o1 )}\leq C.
\ee

The inequality (6) and the standard elliptic estimates yield that
$u_{\e1}-\bar{u}_{\e1}\to G(x,p)$ in $C^{\infty}(\o1 )$ for any
$\o1\subset\subset M\setminus\{p\}$. This completes the proof
of Lemma 2.6.

\hfill $\square$

For any $R>0$, we set $r_{\e1}=\frac{R}{\sl}$.

\begin{lemma}
In $M\setminus B_{r_{\e1}}(0)$, we have
$$
u_{\e1}\geq
G-\xl -2\log (\frac{1+\pi h(p)R^2}{R^2}) - A(p)  + o_{\e1}(1)
$$
where $o_{\e1}(1)\to 0$ as $\e1 \to 0$.
\end{lemma}

Proof:
It is clear that we have
$\btu (u_{\e1}-G-C_{\e1})\leq 0$ for any constant $C_{\e1}$. We choose
$C_{\e1}$ such that
$$
(G+C_{\e1})|_{\partial B_{r_{\e1}}}\leq u_{\e1}|_{\partial B_{r_{\e1}}}.
$$
By Lemma 2.5 and (3) we get
$$
C_{\e1}=-\xl -2\log (\frac{1+\pi h(p)R^2}{R^2}) - A(p) + o_{\e1}(1).
$$
Then the lemma follows from the maximum principle.

\hfill $\square$

Now we are ready to finish the proof of Lemma 2.2.

{\it Proof of Lemma 2.2}: We let $\d1 >0$ small enough so that (3) holds in
$B_{\d1}(p)$. We denote by $o_{\e1}(1)$ (resp. $o_{R}(1)$;
$o_{\d1}(1)$) the terms which tend to 0 as $\e1 \to 0$
(resp. $R\to \infty$; $\d1\to 0$).

We recall that $r_{\e1}=\frac{R}{\sl}$ ($R>0$). We assume that
$\e1 $ is so small that
$\d1 > r_{\e1}$. We have
$$
\int _{M}|\btd u_{\e1}|^2 =
\int _{M\setminus B_{\d1}(p)}|\btd u_{\e1}|^2
+ \int _{B_{\d1}(p)\setminus B_{r_{\e1}}(p) }|\btd u_{\e1}|^2
+ \int _{B_{r_{\e1}}(p) }|\btd u_{\e1}|^2.
$$
It is clear that
\ba
\int _{M\setminus B_{\d1}(p)}|\btd u_{\e1}|^{2}
&=& \int _{M\setminus B_{\d1}(p)}|\btd (u_{\e1}-\bar{u}_{\e1})|^2 \no \\
& = &  \int _{M\setminus B_{\d1}(p)}|\btd G|^2 + o_{\e1}(1) \no \\
& = & -\int _{\partial B_{\d1}}G\cdot \frac{\partial G}{\partial n}
+ o_{\e1}(1)+o_{\d1}(1)
\ea
and
\ba
\int _{B_{r_{\e1}}(p) }|\btd u_{\e1}|^{2}
& = &  \int _{B_R(0)}|\btd \0v |^2 + o_{\e1 }(1)\no\\
& = & 16\pi\log (1+\pi h(p)R^2) -16\pi + o_{\e1 }(1) + o_{R}(1),
\ea
by Lemma 2.6 and Lemma 2.5.

It remains to estimate
$\int _{B_{\d1}(p)\setminus B_{r_{\e1}}(p) }|\btd u_{\e1}|^2$.

Since $u_{\e1}$ satisfies (5),
we have
\bs
\int _{B_{\d1}(p)\setminus B_{r_{\e1}}(p) }|\btd u_{\e1}|^2
& = & -(8\pi -\e1 )\int _{B_{\d1}(p)\setminus B_{r_{\e1}}(p) }u_{\e1}\\
&  & +(8\pi -\e1 )
\int _{B_{\d1}(p)\setminus B_{r_{\e1}}(p) }he^{u_{\e1}} u_{\e1}\\
&  & +\int _{\partial B_{\d1}(p)}u_{\e1}\frac{\partial u_{\e1}}{\partial n}
-\int _{\partial B_{r_{\e1}}(p)}u_{\e1}\frac{\partial u_{\e1}}{\partial n}.
\es
Using Lemma 2.9, we have

\bs
\int _{B_{\d1}(p)\setminus B_{r_{\e1}}(p) }he^{u_{\e1}} u_{\e1}
& \geq & - \xl
\int _{B_{\d1}(p)\setminus B_{r_{\e1}}(p) }he^{u_{\e1}} \\
&  & + \int _{B_{\d1}(p)\setminus B_{r_{\e1}}(p) }he^{u_{\e1}} G
+o_{\e1}(1)+o_{\d1}(1).
\es
Using the equation (5) and the Green formula, one gets that
\bs
(8\pi -\e1 )\int _{B_{\d1}(p)\setminus B_{r_{\e1}}(p) }he^{u_{\e1}} G
& = & -8\pi \int _{B_{\d1}(p)\setminus B_{r_{\e1}}(p) }u_{\e1}
      -\int _{\partial B_{\d1}(p)}\frac{\partial u_{\e1}}{\partial n}G \\
&  & +\int _{\partial B_{\d1}(p)}u_{\e1}\frac{\partial G }{\partial n}
   +\int _{\partial B_{r_{\e1}}(p)}\frac{\partial u_{\e1}}{\partial n} G \\
& & -\int _{\partial B_{r_{\e1}}(p)}u_{\e1}\frac{\partial G}{\partial n}
     +(8\pi -\e1 )\int _{B_{\d1}(p)\setminus B_{r_{\e1}}(p) } G .
\es
By Lemma 2.6 we have
\bs
(8\pi -\e1 )\int _{B_{\d1}(p)\setminus B_{r_{\e1}}(p) }he^{u_{\e1}} G
& = &  -8\pi \int _{B_{\d1}(p)\setminus B_{r_{\e1}}(p) }u_{\e1}
+ \int _{\partial B_{\d1}(p)}\frac{\partial G }{\partial n}\bar{u}_{\e1}\\
& &+\int _{\partial B_{r_{\e1}}(p)}\frac{\partial u_{\e1}}{\partial n} G
 -  \int _{\partial B_{r_{\e1}}(p)}u_{\e1}\frac{\partial G}{\partial n}\\
& & +o_{\e1}(1) + o_{\d1}(1).
\es
Using the equation (5) we also have
\bs
- (8\pi -\e1 )\xl
\int _{B_{\d1}(p)\setminus B_{r_{\e1}}(p) }he^{u_{\e1}}
& = & -(8\pi - \e1 )(Vol (B_{\d1}(p))-Vol (B_{r_{\e1}}(p)))\xl\\
&  & -\xl \int _{\partial B_{r_{\e1}}(p)}\frac{\partial u_{\e1}}{\partial n}
+ \xl \int _{\partial B_{\d1 }(p)}\frac{\partial u_{\e1}}{\partial n}.
\es

We conclude that
\bs
\int _{B_{\d1}(p)\setminus B_{r_{\e1}}(p) }|\btd u_{\e1}|^2
&\geq & -(16\pi -\e1 )\int _{B_{\d1}(p)\setminus B_{r_{\e1}}(p) }u_{\e1}
+ \int _{\partial B_{\d1}(p)}u_{\e1}\frac{\partial u_{\e1}}{\partial n}\\
&  &- \int _{\partial B_{r_{\e1}}(p)}u_{\e1}\frac{\partial u_{\e1}}{\partial n}
+ \int _{\partial B_{\d1}(p)}\frac{\partial G }{\partial n}\bar{u}_{\e1}\\
&  &+\int _{\partial B_{r_{\e1}}(p)}\frac{\partial u_{\e1}}{\partial n} G
 -  \int _{\partial B_{r_{\e1}}(p)}u_{\e1}\frac{\partial G}{\partial n}\\
&  &- \xl \int _{\partial B_{r_{\e1}}(p)}\frac{\partial u_{\e1}}{\partial n}
+ \xl \int _{\partial B_{\d1 }(p)}\frac{\partial u_{\e1}}{\partial n}\\
&  & -(8\pi - \e1 )(Vol (B_{\d1}(p))-Vol (B_{r_{\e1}}(p)))\xl\\
&  &+ o_{\e1}(1) + o_{\d1}(1).
\es
Applying Lemma 2.5 and Lemma 2.9 one has
\bs
\lefteqn{-\int _{\partial B_{r_{\e1}}(p)}\frac{\partial u_{\e1}}{\partial n}
(u_{\e1} -(G-\xl ))}\\
& \geq & \frac{8\pi ^{2}h(p)R^2}{(1+\pi h(p)R^2)}(-A(p)
 -2\log (\frac{1+\pi h(p)R^2}{R^2}))\\
& & +o_{\e1 }(1)+o_R(1).
\es
Using Lemma 2.5 we have
\bs
-\int _{\partial B_{r_{\e1}}(p)}u_{\e1}\frac{\partial G}{\partial n}
& = & -\xl \int _{\partial B_{r_{\e1}}(p)}\frac{\partial G}{\partial n}
- 16\pi\log (1+\pi h(p)R^2)\\
& &  +o_{\e1 }(1)+o_R(1)\\
& = & 8\pi (1- Vol (B_{r_{\e1}}(p)))\xl
- 16\pi\log (1+\pi h(p)R^2)\\
&  & + o_{\e1 }(1)+o_R(1).
\es
By the equation (5) one gets
\bs
\xl \int _{\partial B_{\d1 }(p)}\frac{\partial u_{\e1}}{\partial n}
& = & -(8\pi -\e1 )(1- Vol (B_{\d1}(p)))\xl\\
&  & +(8\pi -\e1 )\xl
\int _{M\setminus B_{\d1}(p)}he^{u_{\e1}} \\
&  &\geq -(8\pi -\e1 )(1- Vol (B_{\d1}(p)))\xl .
\es
Similarly
\bs
\bar{u}_{\e1} \int _{\partial B_{\d1 }(p)}\frac{\partial u_{\e1}}{\partial n}
& = & -(8\pi -\e1 )(1- Vol (B_{\d1}(p)))\bar{u}_{\e1}\\
&  &+ (8\pi -\e1 )\bar{u}_{\e1}e^{\bar{u}_{\e1}}
\int _{M\setminus B_{\d1}(p)}he^{u_{\e1}-\bar{u}_{\e1}} \\
& = & -(8\pi -\e1 )(1- Vol (B_{\d1}(p)))\bar{u}_{\e1} +o_{\e1}(1)
\es
and
$$
\int _{\partial B_{\d1}(p)}\frac{\partial G }{\partial n}\bar{u}_{\e1}
= -8\pi (1- Vol (B_{\d1}(p)))\bar{u}_{\e1}.
$$

We have
\bs
\int _{B_{\d1}(p)\setminus B_{r_{\e1}}(p) }u_{\e1}
&=&
\int _{B_{\d1}(p)\setminus B_{r_{\e1}}(p) }
(u_{\e1}-\bar{u}_{\e1})\\
&  &+(Vol(B_{\d1}(p))-Vol(B_{r_{\e1}}(p)))\bar{u}_{\e1}.
\es
By Lemma 2.3 we have
$$
\int _{B_{\d1}(p)\setminus B_{r_{\e1}}(p) }u_{\e1}
=(Vol(B_{\d1}(p))-Vol(B_{r_{\e1}}(p)))\bar{u}_{\e1}
+o_{\d1}(1).
$$
It is clear that
$$
\xl Vol(B_{r_{\e1}}(p))=o_{\e1}(1).
$$
By Lemma 2.9 we also have
$$
-\bar{u}_{\e1} Vol(B_{r_{\e1}}(p))=o_{\e1}(1).
$$
We therefore have
\ba
\int _{B_{\d1}(p)\setminus B_{r_{\e1}}(p) }|\btd u_{\e1}|^2
& \geq & \e1 \xl -(16\pi -\e1 )\bar{u}_{\e1}
-16\pi\log(1+\pi h(p)R^2)\no\\
&  & +\int _{\partial B_{\d1}}G\cdot \frac{\partial G}{\partial n}
-8\pi A(p) - 16\pi\log \pi -16\pi\log h(p) \no\\
& & +o_{\e1}(1) +o_{R}(1)+o_{\d1}(1).
\ea

It follows from (7), (8) and (9) that
\bs
\int _{M}|\btd u_{\e1}|^2
&\geq &  \e1 \xl -(16\pi -\e1 )\bar{u}_{\e1}\\
&  & -8\pi A(p) - 16\pi\log \pi -16\pi -16\pi\log h(p)\\
& & +o_{\e1}(1) +o_{R}(1)+o_{\d1}(1).
\es
So,
\bs
J_{\e1}(u_{\e1}) & \geq &
\frac{\e1}{2} \xl - \frac{\e1 }{2}\bar{u}_{\e1}
- 4\pi A(p) - 8\pi\log \pi \\
&  &- 8\pi -8\pi\log h(p)
+o_{\e1}(1) +o_{R}(1)+o_{\d1}(1).
\es
Thus, we have
$$
J_{\e1}(u_{\e1})  \geq
- 4\pi A(p) - 8\pi\log \pi
 -  8\pi -8\pi\log h(p)
+o_{\e1}(1) +o_{R}(1)+o_{\d1}(1).
$$
Hence
\be
\inf_{\e1 >0}\l1 _{\e1}
\geq -8\pi -8\pi\log\pi - 4\pi
( \max_{p\in M} (A(p)+2\log h(p))).
\ee
The lemma follows.

\hfill $\square$

Consequently, we have the following lemma which will be used
in the proof of Theorem 1.2.

\begin{lemma}
Assume that the minimizing sequence $u_{\e1}$ of $J_{\e1}$
does not converge in $L_1^2(M)$. Then
$$
\inf _{H_1} J(u)\geq -8\pi -8\pi\log\pi - 4\pi
( \max_{p\in M} (A(p)+2\log h(p))).
$$
\end{lemma}

Proof:
Otherwise, there would exist $u\in H_1$ and $\gamma >0$
such that
$$
J(u)< -8\pi -8\pi\log\pi - 4\pi
( \max_{p\in M} (A(p)+2\log h(p))) -2\gamma .
$$
So,
$$
J_{\e1}(u)< -8\pi -8\pi\log\pi - 4\pi
( \max_{p\in M} (A(p)+2\log h(p))) -\gamma,
$$
when $\e1$ is sufficiently small, which contradicts (10).

\hfill $\square$

\section{Existence theorems}

One can directly prove the following theorem using Lemma 2.10,
because $J(0)=-8\pi\log\int_Mh$.

\begin{theorem}
Let $M$ be a compact Riemann surface.
Let $h(x)$ be a positive smooth function on $M$.
Suppose that
$$
\log\int_Mh > (1+\log\pi )+\frac{1}{2}
\max_{p\in M}(A(p)+2\log h(p)).
$$
Then the equation (1) has a smooth solution.
\end{theorem}

\begin{remark}
If $h$ is a positive constant, then the condition of Theorem 3.1
is satisfied precisely if
$$
\max A(p) < - 2 -2\log \pi.
$$
If $M$ is the standard sphere with volume 1, the constant $A$
in the local expression of G $($see $(3)$$)$ is $-2 -2\log \pi$, and so the
preceding inequality does not hold. We shall see in Section 4, that
it holds for some, but not for all flat tori with volume 1.
\end{remark}

In the sequel, we shall use

\begin{proposition}
Let $M$ be a compact Riemann surface.  Let $K(p)$ be the Gauss
curvature of $M$ at $p$. Let $G(x,p)$ be the
Green function on $M$ satisfying $(2)$. Let $G$ be locally
expressed by $(3)$. Then
$$
c_1+c_3+\frac{2}{3}K(p)=4\pi.
$$
\end{proposition}

Proof:
We denote by $(r,\t1 )$ the chosen normal coordinate system around
$p$. We write $ds^2=dr^2+g^2(r,\t1 )d\t1 ^2$. It is well-known that
$$
g(r,\t1 )=r-\frac{K(p)}{6}r^3+O(r^4).
$$
By the divergence theorem, we have
$$
\int _{\p1 B_{r}}\frac{\p1 G}{\p1 n}
= -8\pi (1-Vol (B_{r})).
$$
So,
$$\int _0^{2\pi}(-\f1{4}{r} + 2c_1 r\cos^2\t1 + 2c_3 r\sin^2\t1 )
(r-\f1{1}{6}K(p)r^3 + O(r^4))d\t1
$$
$$
= -8\pi (1-\pi r^2 + O(r^3)).
$$
Comparing the coefficients of $r^2$, we get
$$
c_1+c_3+\frac{2}{3}K(p)=4\pi.
$$
This proves the proposition.

\hfill $\square$

We now turn to finish the proof of our main theorem, Theorem 1.2.

{\it Proof of Theorem 1.2}:
We shall construct a blow up sequence $\ph$ with
$$
J(\ph )<-8\pi -8\pi\log\pi - 4\pi
( \max_{p\in M} (A(p)+2\log h(p)))
$$
for $\e1$ sufficiently small.
Note that $J(u)=J(u+C)$  for any constant $C$. Combining the
above fact and Lemma 2.10 one gets Theorem 1.2.
Therefore, it only remains to construct the blow up sequence.

Suppose that $A(p)+2\log h(p) =  \max_{x\in M} (A(x)+2\log h(x))$.
Let $r=dist (x,p)$.
We set
$$
\om = -2\log (r^2 +\e1 ),
$$
$$
G=-4\log r + A(p) + b_1r\cos\t1 + b_2r\sin\t1 + \b1,
$$
where $b_1$ and $b_2$ are constants in (3).
$$
\ph =\left\{ \begin{array}{ll}
 \om + b_1r\cos\t1 + b_2r\sin\t1 + \log \e1,  & \quad r\leq \a1 \s1, \\
 (G-\eta \b1 )+C_{\e1}+\log\e1,  & \quad  \a1 \s1 \leq r\leq 2\a1 \s1,\\
 G + C_{\e1 }+\log \e1, & \quad   r\geq 2\a1 \s1.\\
\end{array}
\right.
$$
Here $\eta\in C^{\infty}_0(B_{2\a1\s1 }(p))$ is a cutoff function,
$\eta =1$ in $B_{\a1\s1 }(p)$, $|\btd \eta |\leq \frac{C}{\a1\s1}$,
$$
C_{\e1}=-2\log(\f1{\a1^2 +1}{\a1^2})-A(p)
$$
and $\a1=\a1 (\e1 )$ will be fixed later on satisfying $\a1\to \infty$
and $\a1^2\e1\to 0$ as $\e1\to 0$.

By a simple calculation one has
\bs
\int _{B_{\a1\s1 }}|\btd \ph |^2
&=& 16\pi\log(\a1^2 +1) - 16\pi\f1{\a1^2}{\a1^2 +1}\\
&& -\f1{16\pi}{6}K(p)\a1^2\e1 + \f1{32\pi}{6}K(p)\e1\log(\a1^2 +1)\\
&&+ \pi (b_1^2 + b_2^2 )\a1^2\e1 + O(\e1 ) + O(\a1^4\e1^2).
\es
It is clear that
\bs
\int _{M\setminus B_{\a1\s1 }}|\btd \ph |^2
& = & \int _{M\setminus B_{\a1\s1 }}|\btd G|^2
+ \int _{B_{2\a1\s1 }\setminus B_{\a1\s1 }}|\btd (\eta\b1 )|^2\\
&  & -2 \int _{B_{2\a1\s1 }\setminus B_{\a1\s1 }}
\btd G\cdot \btd(\eta\b1 )\\
& = & -\int _{\p1 B_{\a1\s1 }}G\cdot \f1{\p1 G}{\p1 n}
-8\pi \int _{M\setminus B_{\a1\s1 }} G\\
&  &+ 2\int _{\p1 B_{\a1\s1 }} \f1{\p1 G}{\p1 n}(\eta\b1 )\\
&  & +O(\a1^4\e1^2 ).
\es

Using (3) one has locally
\bs
G(r,\t1 ) &=& -4\log r + A(p) + b_1r\cos\t1 + b_2r\sin\t1 \\
&&+ c_1 r^2\cos^2\t1
+2c_2r^2\cos\t1 \sin\t1 + c_3r^2\sin^2\t1 +O(r^3)
\es
and
\bs
\f1{\p1 G}{\p1 r} &=&
\f1{-4}{r} + b_1\cos\t1 + b_2\sin\t1 \\
&& +2c_1 r\cos^2\t1
+4c_2 r\cos\t1 \sin\t1 + 2c_3r\sin^2\t1 +O(r^2).
\es
So,
\bs
-\int _{\p1 B_{\a1\s1 }}G\cdot \f1{\p1 G}{\p1 n}
& = & -16\pi\log(\a1^2\e1 ) + 8\pi A(p) +
4\pi (c_1 + c_3 )\a1 ^2\e1\\
&  &+ 4\pi (c_1 + c_3 )\a1 ^2\e1\log(\a1 ^2\e1 )
- 2\pi A(p) (c_1 + c_3 )\a1 ^2\e1\\
&  &+ \f1{16\pi}{6} K(p)\a1 ^2\e1\log(\a1 ^2\e1 )
- \f1{8\pi}{6} K(p) A(p) \a1 ^2\e1\\
&  &- \pi (b_1^2 + b_2^2 )\a1^2\e1 + O(\a1 ^4\e1^2\log (\a1 ^2\e1 )).
\es
Similarly,
$$
2\int _{\p1 B_{\a1\s1 }}\eta\b1 \cdot \f1{\p1 G}{\p1 n}
= -8\pi (c_1 + c_3 )\a1 ^2\e1 + O(\a1 ^4\e1 ^2 ).
$$
Hence,
\bs
\lefteqn{-\int _{\p1 B_{\a1\s1 }}G\cdot \f1{\p1 G}{\p1 n}
+2\int _{\p1 B_{\a1\s1 }}\eta\b1 \cdot \f1{\p1 G}{\p1 n}}\\
& = & -16\pi\log(\a1^2\e1 ) + 8\pi A(p) +
4\pi (c_1 + c_3 )\a1 ^2\e1\\
& &+ 4\pi (c_1 + c_3 )\a1 ^2\e1\log(\a1 ^2\e1 )
- 2\pi A(p) (c_1 + c_3 )\a1 ^2\e1\\
&  & +\f1{16\pi}{6} K(p) \a1 ^2\e1\log(\a1 ^2\e1 )
- \f1{8\pi}{6} K(p) A(p) \a1 ^2\e1\\
&  &- 8\pi (c_1 + c_3 )\a1 ^2\e1 - \pi (b_1^2 + b_2^2 )\a1^2\e1
+O(\a1 ^4\e1^2\log (\a1 ^2\e1 )).
\es
Since $\int_M G =0$, we have
\bs
-8\pi\int _{M\setminus B_{\a1\s1 }}G
& = & 8\pi \int _{B_{\a1\s1 }}G\\
& = & -16\pi^2\a1^2\e1\log (\a1 ^2\e1 ) + 16\pi^2 \a1 ^2\e1\\
& &+ 8\pi A(p)\a1 ^2\e1 +O(\a1 ^4\e1^2\log (\a1 ^2\e1 )).
\es
So,
\bs
\int _{M}|\btd \ph |^2
& = & 16\pi\log (\f1{\a1^2+1}{\a1^2}) - 16\pi\log\e1 \\
&  & -16\pi + \f1{16\pi }{\a1^2+1} + 8\pi A(p) \\
&  &+ 4\pi ((c_1 + c_3 ) + \f1{2}{3}K(p) -4\pi )\a1^2\e1\log (\a1 ^2\e1 )\\
&  &+ 4\pi (4\pi -(c_1 + c_3 ) - \f1{2}{3}K(p))\a1 ^2\e1\\
&  &+ 2\pi A(p) (4\pi - (c_1 + c_3 ) -\frac{2}{3}K(p))\a1 ^2\e1 \\
&  &+ \f1{32\pi}{6}K(p)\e1 \log(\a1^2 + 1) +O(\a1 ^4\e1^2\log (\a1 ^2\e1 )).
\es
Applying Proposition 3.3, one has
\bs
\int _{M}|\btd \ph |^2
& = & 16\pi\log (\f1{\a1^2+1}{\a1^2}) - 16\pi\log\e1 \\
&  &- 16\pi + \f1{16\pi }{\a1^2+1} + 8\pi A(p) \\
&  &+ \f1{32\pi}{6} K(p) \e1 \log(\a1^2 + 1)  + O(\a1 ^4\e1^2\log (\a1 ^2\e1
)).
\es
Calculating directly, one has

\bs
\int _{B_{\a1\s1 }}\om
& = &-2\pi \a1 ^2\e1\log (\a1^2+1)\e1 - 2\pi\e1\log (\a1^2+1)\\
&  & +2\pi \a1 ^2\e1 + O(\a1 ^4\e1^2\log (\a1 ^2\e1 )).
\es
It is also obvious that
\bs
\int _{M\setminus B_{\a1\s1 }}\ph
& = & (1- Vol (B_{\a1\s1 }))\log\e1 - \int _{B_{\a1\s1 }}G\\
&  &+ C_{\e1} (1- Vol (B_{\a1\s1 })) - \int _{B_{2\a1\s1 }\setminus B_{\a1\s1}}
\eta\b1 \\
& = & 2\pi \a1 ^2\e1\log (\a1^2\e1) - 2\pi \a1^2\e1 - A(p)\pi \a1^2\e1\\
&  &+ (1- Vol (B_{\a1\s1 }))\log\e1 +  C_{\e1} (1- Vol (B_{\a1\s1 }))
+ O(\a1^4\e1^2 ).
\es
Thus,
\bs
\int _{M}\ph & = & \log\e1 -2\pi \a1 ^2\e1\log (\frac{\a1^2+1}{\a1^2})
-2\pi \e1\log (\a1 ^2 +1)\\
&  &- A(p) -2 \log (\frac{\a1^2+1}{\a1^2})(1- Vol (B_{\a1\s1 }))\\
&  &+ O(\a1 ^4\e1^2\log (\a1 ^2\e1 )).
\es
We have
\bs
\int _{B_{\a1\s1 }}e^{\ph}
& = & \e1 \int _{0}^{\a1\s1 }\frac{2\pi}{(r^2+\e1 )^2}(r-\f1{1}{6} K(p)r^3)dr\\
&  &+ \f1{\e1}{2}\int _0^{\a1\s1 }\int _0^{2\pi}\f1{b_1^2r^2\cos^2\t1 +
b_2^2r^2\sin^2\t1 }{(r^2 +\e1 )^2}d\t1 rdr + O(\e1 )\\
& = & \pi (\f1{\a1 ^2}{\a1 ^2 +1}) - \f1{1}{6}K(p)\pi\e1\log (\a1^2+1) \\
&  & +\f1{\pi}{4}(b^2_1 + b_2^2 ) \e1 \log (\a1 ^2 +1 ) + O(\e1 ).
\es

We choose $\d1 >0$ sufficiently small so that $G$ has the expression
(3) in $B_{\d1}(p)$, then we have

\bs
\int _{M\setminus B_{\a1\s1 }}e^{\ph}
& = & \e1 \int _{M\setminus B_{\d1}}e^{G+C_{\e1}}
+\e1 \int _{B_{\d1}\setminus B_{\a1\s1}}e^{-4\log r + A(p) + C_{\e1}}\\
&  & +\e1 \int _{B_{2\a1\s1 }\setminus B_{\a1\s1 }}e^{G+C_{\e1}}
(e^{-\eta\b1}-1 )\\
&  & +\e1 \int _{B_{\d1}\setminus B_{\a1\s1}}e^{-4\log r + A(p) + C_{\e1}}
(e^{ b_1r\cos\t1 + b_2r\sin\t1 + \b1}-1).
\es
Calculating directly one gets
\bs
\e1 \int _{B_{\d1}\setminus B_{\a1\s1}}e^{-4\log r + A(p) + C_{\e1}}
& =  & \pi \f1{\a1^2}{(\a1^2 +1)^2} \\
& &+ 2\pi \f1{\a1^4}{(\a1^2+1)^2}\f1{1}{6}K(p)
\e1\log (\a1\s1 ) + O(\e1 )
\es
and
\bs
\lefteqn{\e1 \int _{B_{\d1}\setminus B_{\a1\s1}}e^{-4\log r + A(p) + C_{\e1}}
(e^{ b_1r\cos\t1 + b_2r\sin\t1 + \b1}-1)}\\
& = &-\pi\f1{\a1^4}{(\a1^2 + 1)^2}(c_1 + c_3 )\e1 \log (\a1\s1)\\
&  & -\f1{\pi}{2}\f1{\a1^4}{(\a1^2+1)^2}(b_1^2+b_2^2)\e1\log (\a1\s1 ) + O(\e1
).
\es

We therefore have
\bs
\int _{M\setminus B_{\a1\s1 }}e^{\ph}
& = & \pi \f1{\a1^2}{(\a1^2 +1)^2}
 + 2\pi \f1{\a1^4}{(\a1^2+1)^2}\f1{1}{6}K(p)
\e1\log (\a1\s1 )\\
&  & -\pi\f1{\a1^4}{(\a1^2 + 1)^2}(c_1 + c_3 )\e1 \log (\a1\s1)\\
& & -\f1{\pi}{2}\f1{\a1^4}{(\a1^2+1)^2}(b_1^2+b_2^2)\e1\log (\a1\s1 ) + O(\e1
).
\es
Thus,
\bs
\int _{M}e^{\ph}
& = & \pi \f1{\a1^2}{\a1^2 +1} (1+ \f1{1}{\a1^2 +1}
-\f1{\a1^2+1}{\a1^2}\f1{1}{6}K(p)
\e1\log (\a1^2 +1)\\
&  &- \f1{\a1^2}{\a1^2 + 1}(c_1 + c_3 - \frac{1}{3}K(p)) \e1 \log (\a1\s1)\\
&  &+ \f1{1}{4}\f1{\a1^2+1}{\a1^2}(b^2_1 + b_2^2 ) \e1 \log (\a1 ^2 +1 )\\
&  & -\f1{1}{2}\f1{\a1^2}{\a1^2+1}(b_1^2+b_2^2)\e1\log (\a1\s1 )) + O(\e1 ).
\es
It is clear that
$$
\int _M he^{\ph } = h(p)\int _M e^{\ph } +
\int _M (h - h(p))e^{\ph }.
$$

Suppose that
\bs
h(x)-h(p) & = & k_1r\cos\t1 + k_2r\sin\t1 \\
& + & k_3r^2\cos^2\t1 + 2k_4r^2\cos\t1 \sin\t1
+ k_5r^2\sin^2\t1 + O(r^3 )
\es
in $B_{\d1}(p)$.

By a simple computation, we obtain
\bs
\int _{B_{\a1\s1}}(h - h(p))e^{\ph }
& = & \f1{\pi}{2}(k_3+k_5)\e1\log (\a1 ^2 + 1)\\
&  & +\f1{\pi}{2}(k_1b_1 + k_2b_2)\e1\log (\a1 ^2 + 1) + O(\e1 )
\es
and
\bs
\int _{M\setminus B_{\a1\s1}}(h - h(p))e^{\ph }
& = & \int _{B_{\d1}\setminus B_{\a1\s1}}(h - h(p))e^{\ph }
+ \int _{M\setminus B_{\d1}}(h - h(p))e^{\ph }\\
& = & -\f1{\pi}{2}(k_3+k_5)(\f1{\a1^2}{\a1^2+1})^2
\e1\log (\a1^2\e1 ) \\
&   &- \f1{\pi}{2}(k_1b_1 + k_2b_2)(\f1{\a1^2}{\a1^2+1})^2
\e1\log (\a1^2\e1 ) + O(\e1 ).
\es
So,
\bs
\int _{M}(h - h(p))e^{\ph }
& = & \f1{\pi}{4}(\btu h(p))\e1\log (\a1 ^2 + 1)\\
&  & +\f1{\pi}{2}(k_1b_1 + k_2b_2)\e1\log (\a1 ^2 + 1)\\
&  & -\f1{\pi}{4}(\btu h(p))(\f1{\a1^2}{\a1^2+1})^2\e1\log (\a1^2\e1 )\\
&  & -\f1{\pi}{2}(k_1b_1 + k_2b_2)(\f1{\a1^2}{\a1^2+1})^2
\e1\log (\a1^2\e1 ) + O(\e1 ).
\es
Therefore,
\bs
\int _M he^{\ph }
\!\!\!\!& = & h(p) \pi \f1{\a1^2}{\a1^2 +1} (1+ \f1{1}{\a1^2 +1} -
\f1{\a1^2+1}{\a1^2}\f1{1}{6}K(p)\e1\log (\a1^2 +1)\\
&  &+ \f1{1}{4}\f1{\a1^2+1}{\a1^2}(b^2_1 + b_2^2 ) \e1 \log (\a1 ^2 +1 )
 -  \f1{1}{2}\f1{\a1^2}{\a1^2+1}(b_1^2+b_2^2)\e1\log (\a1\s1 ) \\
&  & -\f1{\a1^2}{\a1^2 + 1}(c_1 + c_3 - \frac{1}{3}K(p) )\e1 \log (\a1\s1))
+  \f1{\pi}{4}(\btu h(p))\e1\log (\a1 ^2 + 1)\\
&  & -\f1{\pi}{4}(\btu h(p))(\f1{\a1^2}{\a1^2+1})^2
\e1\log (\a1^2\e1 ) + \f1{\pi}{2}(k_1b_1 + k_2b_2)\e1\log (\a1 ^2 + 1)\\
&  & -\f1{\pi}{2}(k_1b_1 + k_2b_2)(\f1{\a1^2}{\a1^2+1})^2
\e1\log (\a1^2\e1 )  + O(\a1^4\e1^2 )+ O(\e1 ).
\es

Adding the terms in the functional, we get
\bs
J(\ph )
\!\!\!\!& = & -8\pi - 8\pi\log \pi - 4\pi A(p) - 8\pi\log h(p)\\
& &- 16\pi^2 (1-\f1{1}{4\pi}K(p))\e1\log (\a1^2+1)
+ 4\pi (c_1 + c_3 - \frac{1}{3}K(p) )\e1 \log (\a1^2\e1 )\\
& &- 2\pi (b^2_1 + b_2^2 ) \e1 \log (\a1 ^2 +1 )
+  2\pi (b_1^2+b_2^2)\e1\log (\a1^2\e1 )\\
& &- 4\pi \f1{(k_1b_1 + k_2b_2)}{h(p)} \e1 \log (\a1 ^2 +1 )
+  4\pi \f1{(k_1b_1 + k_2b_2)}{h(p)} \e1\log (\a1^2\e1 )\\
& &- 2\pi \f1{\btu h(p)}{h(p)}\e1\log (\a1 ^2 + 1)
+ 2\pi \f1{\btu h(p)}{h(p)}\e1\log (\a1^2\e1 ) \\
& &+ O(\f1{\e1\log (\a1^2+1)}{\a1^2})
+ O(\f1{-\e1 \log (\a1^2\e1 )}{\a1^2})\\
& &+ O((\e1 \log (\a1^2 + 1 ))^2) + O((-\e1 \log (\a1^2\e1 ))^2)\\
& &+ O(\f1{1}{\a1^4}) + O(\a1^4\e1^2 \log (\a1^2\e1 )) + O(\e1 ).
\es
Choosing $\a1$ so that $\a1 ^4\e1 = \f1{1}{\log(-\log\e1 )}$
and applying Proposition 3.3, we get
\bs
J(\ph )
\!\!& = & -8\pi - 8\pi\log \pi - 4\pi A(p) - 8\pi\log h(p)\\
&  &- 16\pi^2 (1-\f1{1}{4\pi} K(p)+\f1{b_1^2+b_2^2}{8\pi}
+\f1{\btu h(p)}{8\pi h(p)} + \f1{(k_1b_1 + k_2b_2)}{4\pi h(p)})\e1(-\log \e1
)\\
&  &+ o(\e1 (-\log \e1)).
\es
This proves theorem 1.2.

\hfill $\square$

\section{ The Green function on a flat torus}

For details on the Green function, we refer to [L].

Let $z=x+iy$ be a variable in $ C$ (the complex plane) and let $\tau=u+iv$,
$v>0$. Here for simplicity, we assume $u=0$. Let $q=e^{-2\pi v}$ and
$q_z=e^{2\pi iz}$. Let $\Sigma_q=\Sigma_v={ C}^*/ (q^{ Z})$,
where $Z$ is the set of integers, ${ C}^*={ C}-\{0\}$
and $q$ acts on ${ C}^*$ by the usual
multiplication. In other words, $\Sigma_q$ is the torus generated by
the lattice $[1, \tau]$. Define a metric on $\Sigma_q$ by
$$ds^2=\frac {1}{v}dx\wedge dy.$$

The area of $\Sigma_q$ with respect to $ds^2$ is $1$.
The corresponding Green function is
$$G(z,0)=-4\log|q^{B_2(y/v)/2}(1-q_z)\prod_{n=1}^{\infty}
(1-q^nq_z)(1-q^{-n}q_z)|,$$
where $B_2(y)=y^2-y+\frac16$ is the second Bernouli polynomial.
Recall the definition of the Green function in the introduction.

Now the asymptotic expansion of the above Green function at the origin is
\bs
\lefteqn{-4\log|z|-4\log 2\pi+\frac{2v\pi}{3}-8\log
(\prod_{n=1}^{\infty}(1-e^{-2\pi nv}))+O(|z|^2)}\\
&=&
-4\log( v^{1/2}|z|) + 2\log v - 4\log 2\pi + \frac{2v\pi}{3}\\
&&- 8\log
(\prod_{n=1}^{\infty}(1-e^{-2\pi nv}))+ O(|z|^2)\\
&=&
-4\log r +2\log v-4\log 2\pi+\frac{2v\pi}{3}\\
&&- 8\log
(\prod_{n=1}^{\infty}(1-e^{-2\pi nv}))+O(|r|^2)
\es
where $r= v^{1/2}|z|$.
The latter expansion is in normal coordinates. Therefore,
$$
A_v=-2\log v-4\log 2\pi+2\frac{v\pi}{3}-8\log
(\prod_{n=1}^{\infty}(1-e^{-2\pi nv}))
$$
Clearly, the asymptotic expansion of the Green function on $\Sigma_v$ is
independent of the base point $0$. $A_v$ is increasing between
$[1,+\infty)$.  Furthermore $A_{1}<-2-2\log \pi=A_0$ and
$\lim_{v\to +\infty}A_v=+\infty$. Hence there exists $v^*\in (1, +\infty)$
 such that

(i) $A_v<A_0$, if $1\le v <v^*,$

(ii) $A_v>A_0$, if $v^*<v<\infty.$

\vspace{.2in}

\begin{center}
{\large\bf REFERENCES}
\end{center}
\footnotesize
\begin{description}

 \item[{[A]}] {Aubin, T., Meilleures constantes dans le theoreme
 d'inclusion de Sobolev et un theoreme de Fredholm non lineaire
 pour la transformation conforme de courbure scalaire, J. Funct.
 Anal. 32(1979), 148-174.}

 \item[{[BM]}] {Brezis, H. and Merle, F., Uniform estimates and blow up
 behavior for solutions of $-\btu u = V(x)e^u $ in two dimensions,
 Comm. Partial Diff. Equat. 16(1991), 1223-1253.}

 \item[{[CaY]}] {Caffarelli, L. and Yang, Y.S., Vortex condensation in the
 Chern-Simons Higgs model: an existence theorem, Comm. Math. Phys. 168(1995),
 321-336.}

 \item[{[CY1]}] {Chang, A. S. Y. and Yang, P., Conformal deformation
 of metrics on $S^2$, J. Diff. Geom. 23(1988), 259-296.}

 \item[{[CY2]}] {Chang, A. S. Y. and Yang, P., Prescribing Gaussian
 curvature on $S^2$, Acta Math. 159(1987), 214-259.}

\item[{[CL]}] {Chang, K. C. and Liu, J. Q., On Nirenberg's problem,
 International J. Math. 4(1993) 35-57.}

 \item[{[CD1]}] {Chen, W. X. and Ding, W. Y. Scalar curvature on $S^2$,
 Trans. Amer. Math. Soc. 303(1987), 365-382.}

 \item[{[CD2]}] {Chen, W. X. and Ding, W. Y., A problem concerning the
 scalar curvature on $S^2$, Kexue Tongbao, Sci. Bull. Ed. 33(1988),
 533-537.}

 \item[{[CL1]}] {Chen, W. X. and Li, C., Prescribing Gaussian
 curvature on surfaces with conical singularities, J. Geom.
 Anal. 1(1991), 359-372.}

 \item[{[CL2]}] {Chen, W. X. and Li, C., Classification of solutions
 of some nonlinear elliptic equations, Duke Math. J. 63(1991),
 615-622.}

 \item[{[D]}] {Ding, W., On the best constant in a Sobolev inequality on
compact 2-manifolds
 and application, unpublished manuscript.}

 \item[{[DJLW]}] {Ding, W., Jost, J., Li, J. and Wang, G., An analysis of the
two-vortex case in the Chern-Simons Higgs model, preprint, June 1997.}

 \item[{[ES]}] {Escobar, J. F. and Schoen, R. M., Conformal metrics
 with prescribed scalar curvature, Invent. Math. 86(1986), 243-253.}

 \item[{[H]}] {Hong, C. W., A best constant and the Gaussian curvature,
 Proc. Amer. Math. Soc. 97(1986), 737-747.}

 \item[{[HKP]}] {Hong, J., Kim, Y. and Pac, P.Y., Multivortex solutions
 of the Abelian Chern-Simons theory, Phys. Rev. Lett. 64(1990), 2230-2233.}

 \item[{[JW]}] {Jackiw, R. and Weinberg, E., Self-dual Chern-Simons vortices,
 Phys. Rev. Lett. 64(1990), 2234-2237.}

 \item[{[KW]}] {Kazdan, J. and Warner, F., Curvature functions for
 compact 2-manifolds, Ann. Math. 99(1974), 14-47.}

 \item[{[L]}] {Lang, S. Introduction to Arakelov theory, Springer-Verlag 1988}

 \item[{[M1]}] {Moser, J., A sharp form of an inequality of N. Trudinger,
 Indiana Univ. Math. J. 20(1971), 1077-1092.}

 \item[{[M2]}] {Moser, J., On a nonlinear problem in differential geometry,
 in: Dynamical Systems, ed. M. Peixoto, Academic Press, 1969, 273-280.}

 \item[{[Sc]}] {Schoen, R. M, Conformal deformation
 of a Riemannian metric to constant scalar curvature, J. Diff. Geom.
 20(1984), 479-495.}

 \item[{[SY]}] {Schoen, R. M and Yau, S. T., On the proof of the positive
 mass conjecture in general relativity, Comm. Math. Phys. 65(1979), 45-76.}

 \item[{[T]}] {Tarantello, G., Multiple condensate solutions for the
 Chern-Simons-Higgs theory, J. Math. Phys. 37(1996), 3769-3796.}

\end{description}
\end{document}